\documentclass{amsart}

\usepackage{amsmath,amssymb,amsthm}
\usepackage{graphicx,tikz}
\newtheorem{theorem}{Theorem}

\theoremstyle{definition}

\theoremstyle{remark}

\begin{document}

\title[]{A Spectral Approach to the \\Shortest Path Problem}
\thanks{S.S. is supported by the NSF (DMS-1763179) and the Alfred P. Sloan Foundation.}

\author[]{Stefan Steinerberger}
\address{Department of Mathematics, Yale University, New Haven, CT 06511, USA}
\email{stefan.steinerberger@yale.edu}

\begin{abstract} 
Let $G=(V,E)$ be a simple, connected graph. One is often interested in a short path between two vertices $u,v$. We propose a spectral algorithm: construct the function $\phi:V \rightarrow \mathbb{R}_{\geq 0}$
$$ \phi = \arg\min_{f:V \rightarrow \mathbb{R} \atop  f(u) = 0, f \not\equiv 0} \frac{\sum_{(w_1, w_2) \in E}{(f(w_1)-f(w_2))^2}}{\sum_{w \in V}{f(w)^2}}.$$
$\phi$ can also be understood as the smallest eigenvector of the Laplacian Matrix $L=D-A$ after the $u-$th row and column have been removed. We start in the point $v$ and construct a path from $v$ to $u$: at each step, we move to the neighbor for which $\phi$ is the smallest. This algorithm provably terminates and results in a short path from $v$ to $u$, often the shortest. The efficiency of this method is due to a discrete analogue of a phenomenon in Partial Differential Equations that is not well understood. We prove optimality for trees and discuss a number of open questions.
\end{abstract}

\vspace{-10pt}

\maketitle

\section{Introduction and Algorithm}

\subsection{Introduction.} 
The Shortest Path Problem asks for the shortest path connecting two vertices in a graph. For simplicity, we will only deal with simple, connected graphs throughout the entire paper. 
The problem is one of the foundational problems in computer science and has received a lot of attention 
\cite{bell, cher, dij, ford, fred0, fred, hart, ing, joh, karl, roy, pet, pol, sei, shi, thor, war, will, whi}.  
 In the case that we are considering (unweighted and undirected), breadth first search is known to cost $\mathcal{O}(|V| + |E|)$. 

\begin{figure}[h!]
\begin{minipage}[l]{.48\textwidth}
\begin{tikzpicture}
\node at (0,0) {\includegraphics[width = 0.8\textwidth]{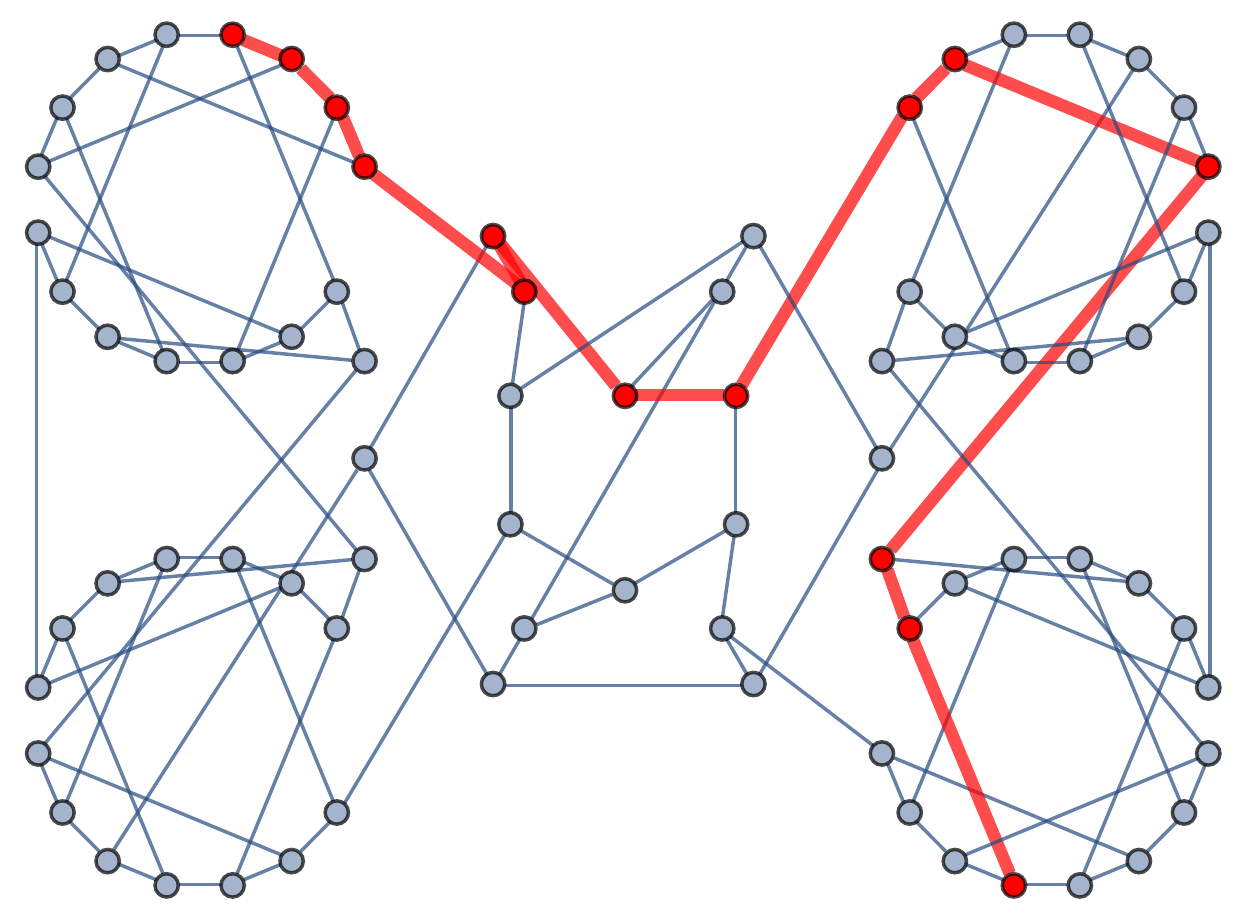}};
\end{tikzpicture}
\end{minipage} 
\begin{minipage}[r]{.48\textwidth}
\begin{tikzpicture}
\node at (0,0) {\includegraphics[width = 0.8\textwidth]{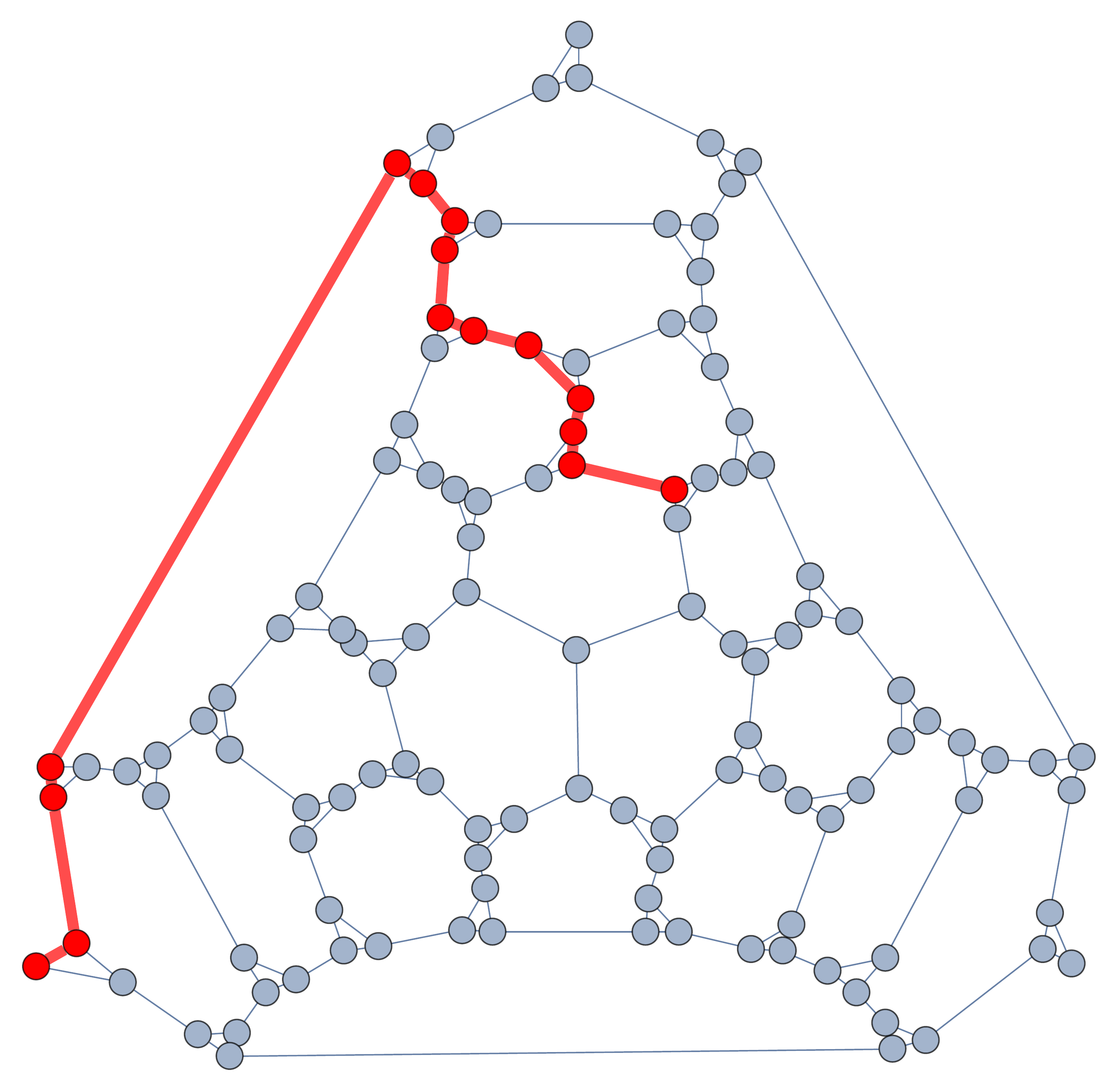}};
\end{tikzpicture}
\end{minipage} 
\caption{Paths taken by the Spectral Method.} 
\end{figure}

We introduce a new algorithm based on Spectral Graph Theory; the main point is not to find a new algorithm to compute the shortest path (though it may have interesting applications, see \S 2.4) but to introduce a phenomenon in Spectral Graph Theory. The algorithm does not always result in the shortest path but it does seem to result in surprisingly short paths (and even on fairly strange graphs it tends to find the shortest path for a \textit{vast} majority of vertices, see \S 3). What is interesting about this algorithm is that the reason for its effectiveness is closely related to some famously unsolved problems in Partial Differential Equations (see \S 2.4): \textit{any} type of result that one can prove about this algorithm on graphs may lead to a better understanding of the PDE problem in the continuous setting! There is an understanding that `robust' phenomena in PDEs also occur on graphs and, conversely, anything that one can prove on graphs one can often also prove for PDEs, see \cite{cald, cald2, trill, calder} for recent examples. One contribution of our paper is to show that the phenomenon that causes `Hot Spots' of eigenfunctions to be on the boundary is `robust' in this sense.
In light of this, any progress on the discrete problem could be of substantial interest.

\subsection{The Algorithm.} We now present the algorithm. We assume that $G=(V,E)$ is a simple, connected graph with $|V| = n$ vertices which, for simplicity, we label $\left\{1, \dots, n \right\}$. We will assume that $D \in \mathbb{R}^{n \times n}$ is a diagonal matrix where $d_{ii} = \mbox{deg}(i)$ is the degree of the vertex $i$. We will use $A$ to denote the (symmetric) adjacency matrix $A = (a_{ij}) \in \left\{0,1 \right\}^{n \times n}$ where $a_{ij} = 1$ if and only if $(i,j) \in E$. 
\begin{quote} \textbf{Algorithm.}\\
Input: a graph $G=(V,E)$ and two vertices $i,j \in V$.\\
Output: a path from $j$ to $i$.

\begin{enumerate}
\item[] (\textit{Preprocessing})
\item Construct the matrix $L=D-A$ and construct $L_i$ by deleting the $i-$th row and the $i-$th column of $L$.
\item Find an eigenvector $\phi \in \mathbb{R}^n$ associated to the smallest positive eigenvalue of $L_i$ and interpret
it as a function $\phi:V\setminus \left\{i\right\} \rightarrow \mathbb{R}$.  
\item If $\phi$ has negative entries, multiply it with $-1$. Set $\phi(i) = 0$.
\item[] (\textit{Path Construction})
\item Construct a path by starting in $x_0 = j$. 
\item While $x_k \neq i$, look at all the neighbors of $x_k$ and pick $x_{k+1}$ to be a neighbor for which $\phi$ assumes its minimum.
\end{enumerate}
\end{quote}

In practical problems the algorithm can be applied twice: once to find a path from $i$ to $j$ and then to find one from $j$ to $i$. We observe in practice that this trick leads to optimal paths in many settings (see \S 3).
 We will now give some simple heuristic motivating the algorithm. For any function $f:V \rightarrow \mathbb{R}$, we can consider the quadratic form 
 $$ \left\langle f, Lf\right\rangle= \sum_{(w_1, w_2) \in E}{(f(w_1)-f(w_2))^2}$$
and try to minimize it. Clearly, the quadratic form has a scaling symmetry under multiplication with constants: it is thus customary to assume 
$$ \sum_{v \in V}{f(v)^2} = 1.$$
For simple, connected graphs, this minimization procedure is quite simple: clearly, the minimum will be given by constant functions resulting in the value 0. This is also sometimes referred to as the `trivial' eigenvector of the Graph Laplacian. 
Instead, what we will do is to pick a fixed $i \in V$ and consider the minimization problem over all functions $f:V \rightarrow \mathbb{R}$ for which $f(i) = 0$, i.e.
$$ \phi = \arg \min_{f(i) = 0} \sum_{(w_1, w_2) \in E}{(f(w_1)-f(w_2))^2} \qquad \mbox{subject to}~ \sum_{v \in V}{f(v)^2} = 1.$$
This minimization problem cannot result in constant functions. The way it is set up, we expect $f$ to be small around vertices close to $i$ and larger when it is further away. This simple heuristic, though meaningful and easily verified in practice, is not really mathematically understood (see \S 2.4.).  In terms of computation, it is easily verified that
$$ f^T (D-A)f =  \sum_{(w_1, w_2) \in E}{(f(w_1)-f(w_2))^2}$$
and thus minimizing the quadratic form is akin to finding the smallest eigenvectors of $D-A$ (which is the constant eigenvector). If we additionally enforce that $f(i) = 0$, then a quick inspection on the left-hand side shows that this the same problem after having removed the $i-$th row and $i-$th column of $D-A$. Since the graph is connected, the smallest eigenvalue will be larger than 0. The smallest eigenvector $\phi$ is unique and does not change sign (so up to multiplication with $-1$ we can assume it to have nonnegative entries). We note that if removing the vertex $i$ disconnects the graph (as is, for example, the case for some vertices on trees), then the problem decouples into two or more problems and can be solved one-by-one on each graph.

\subsection{Numerical Remarks.}
One might a priori assume that minimizing a quadratic form, finding the smallest eigenvalue of a matrix, is prohibitively expensive. However, there have been recent spectacular advances on nearly-in-time solvers of equations involving Graph Laplacians, we refer to work of Spielman \& Teng \cite{spiel1, spiel2, spiel3}, the papers \cite{fast1, fast2,fast3,fast4,fast5,fast6} and the book of Vishnoi \cite{vishnoi}. We point out that Spielman \& Teng \cite{spiel3} give an algorithm that computes an approximation of the first nontrivial eigenfunction in $\widetilde{\mathcal{O}}(|V| + |E|)$.
Similarly expensive techniques have been profitably used in Computer Graphics, see \S 2.4.

\section{The Results}

\subsection{The Algorithm terminates.}
Our first result is that the algorithm provably terminates. This requires us to prove several properties about the solution of the minimization problem (in particular, uniqueness, that it only vanishes in one point and that there are no local minima).
\begin{theorem} The Algorithm produces a path between the vertices. 
\end{theorem}

The argument follows classical arguments from the continuous setting in spirit, we refer to Courant \& Hilbert \cite{cour, cour2}. One downside is that the argument does not seem to result in quantitative estimates which would be desirable: as emphasized above, and in other parts of the paper, this is not merely an interesting question in terms of the algorithm proposed here. A better understanding of the algorithm would likely to lead progress on some notoriously difficult questions in partial differential equations (see \S 2.4).

\subsection{Trees.} We can prove that the algorithm produces the optimal path on graphs that are connected trees.

\begin{theorem} The Algorithm produces the shortest path on connected trees.
\end{theorem}

The argument makes use of a reflection trick that allows us to reduce the problem to one the has been studied previously (with a recent and particularly short solution by Lederman and the author \cite{leder}). Needless to say, it would be very nice to have an extension of the result to broader classes of graphs.

\subsection{A Continuous Analogue.}
In this section, we describe a continuous analogue of the statements above (which, indeed, was motivated by the results above). We will consider a simply connected domain $\Omega \subset \mathbb{R}^2$ with smooth boundary. We assume the boundary is comprised of two parts: a simply connected set $\partial \Omega_0$ (homeomorphic to an interval) and the complement $\partial \Omega_n$ (see Fig. 2 for a sketch of how such a domain could look like).

\begin{center}
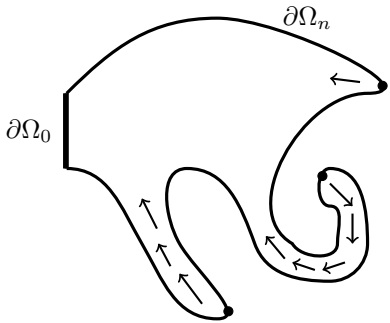
\begin{figure}[h!]
\begin{tikzpicture}
\draw [ultra thick] (0,0) -- (0,1);
\draw [ultra thick] (-0.02,0) -- (-0.02,1);
\draw [very thick] (0,0) to[out=0, in=120] (1,-1) to[out=300, in=180] (2, -2) to[out=0, in=300] (1.5, -1) to[out=120, in=180] (1.6, 0) to[out=0, in=110] (2.5, -1) to[out=290, in =180] (3.5, -1.5) to [out=0, in =270] (4, -1) to[out=90, in=0] (3.5, 0) to[out=180, in = 180] (3.5, -0.5) to[out=0, in=90] (3.6, -1) to[out=270, in=300] (3, -1) to[out=150, in =180] (4, 1) to[out=0, in=0] (2, 2) to[out=180, in=45] (0,1);
\draw [thick, ->]  (1.8,-1.8) -- (1.5,-1.4);
\draw [thick, ->]  (1.38,-1.3) -- (1.22,-1);
\draw [thick, ->]  (1.2,-0.8) -- (1,-0.4);
\filldraw (2.15, -1.9) circle (0.06cm);
\draw [thick, ->] (3.9, 1.15) -- (3.5, 1.2);
\filldraw (4.2, 1.1) circle (0.06cm);
\node at (-0.5, 0.5) {$\partial {\Omega_0}$};
\node at (3.2, 2) {$\partial {\Omega_n}$};
\draw [thick, ->] (3.5, -0.2) -- (3.8, -0.5);
\filldraw (3.4, -0.1) circle (0.06cm);
\draw [thick, ->] (3.8, -0.6) -- (3.8, -1);
\draw [thick, ->] (3.7, -1.25) -- (3.4, -1.35);
\draw [thick, ->] (3.3, -1.35) -- (3, -1.25);
\draw [thick, ->] (2.9, -1.2) -- (2.65, -0.9);
\end{tikzpicture}
\caption{Theorem 3 illustrated: the gradient flow moves towards the boundary and impacts on the part of the boundary where the function vanishes.}
\end{figure}
\end{center}

We will consider the differential inequality
\begin{align*}
\Delta u < 0 \qquad &\mbox{inside} ~\Omega\\
u = 0 \qquad &\mbox{on}~\partial \Omega_0\\
\frac{\partial u}{\partial n} = 0  \qquad &\mbox{on}~\partial \Omega_n,
\end{align*}
where $n$ is the normal vector on the boundary.
This is, in some sense, the canonical analogue of the structure we encounter on graphs: there is a vertex pinned down to 0 on which $\phi$ vanishes (and it vanishes nowhere else), Neumann boundary conditions are the natural analogue to the absence of boundary conditions on the graph. A particularly canonical way of obtaining this continuous setup is to consider the first nontrivial Laplacian eigenfunction $-\Delta u = \mu_2 u$ with Neumann boundary conditions and to take a connected component of $\left\{x:u(x) > 0\right\}$. 
We will now study gradient flows on this space: at each point $x \in \Omega$ where the gradient $\nabla u$ does not vanish, we can define a local vector field pointing in the direction where $u$ decays the fastest 
$$ V(x) = -\frac{\nabla u(x)}{|\nabla u(x)|}.$$
We can now, correspondingly, define a gradient flow started in a point $x_0$ to be the solution of the ODE $\gamma(0) = x_0$ and
$$ \frac{d}{dt} \gamma(t) = V(\gamma(t)).$$
This is the continuous analogue of our way of constructing the path: at each step, we add the adjacent vertex with the smallest value on it.

\begin{theorem} Let $\Omega$ be a simply connected domain and $u$ be as above. If every local maximum of $u$ is on the boundary and the critical points are non-degenerate, then all the critical points are on the boundary. Moreover, the gradient flow started in any point inside $\Omega$ converges to $\partial \Omega_0$. 
\end{theorem}
Very few results in this direction seem to be known, we refer to \cite{alon}.
It could be quite interesting to make this statement quantitative. Is there an upper bound on the length of the gradient flow in terms of geometric quantities of $\Omega$?

\begin{figure}[h!]
\begin{minipage}[l]{.48\textwidth}
\begin{tikzpicture}
\node at (0,0) {\includegraphics[width = 0.85\textwidth]{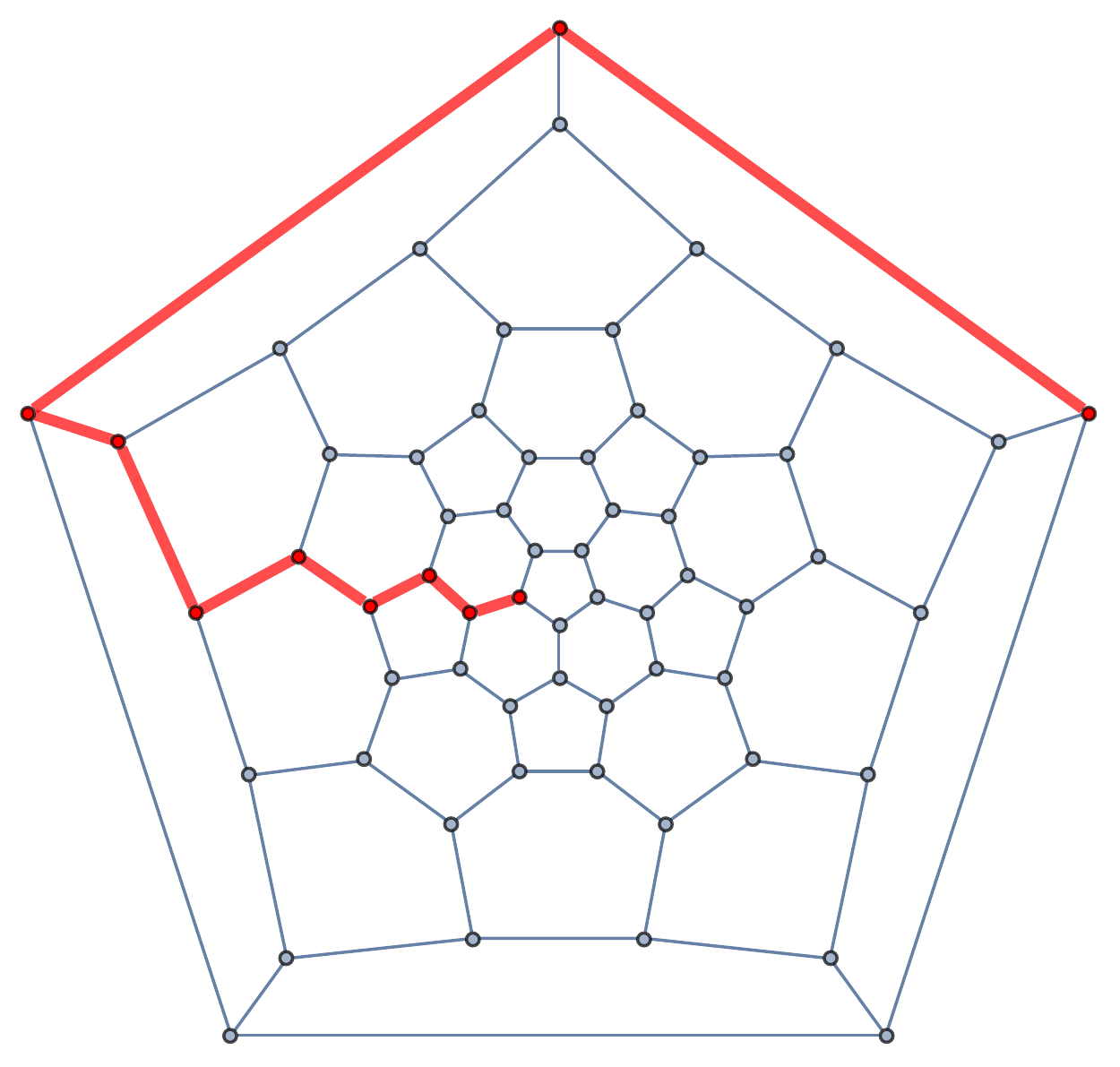}};
\end{tikzpicture}
\end{minipage} 
\begin{minipage}[r]{.48\textwidth}
\begin{tikzpicture}
\node at (0,0) {\includegraphics[width = 0.85\textwidth]{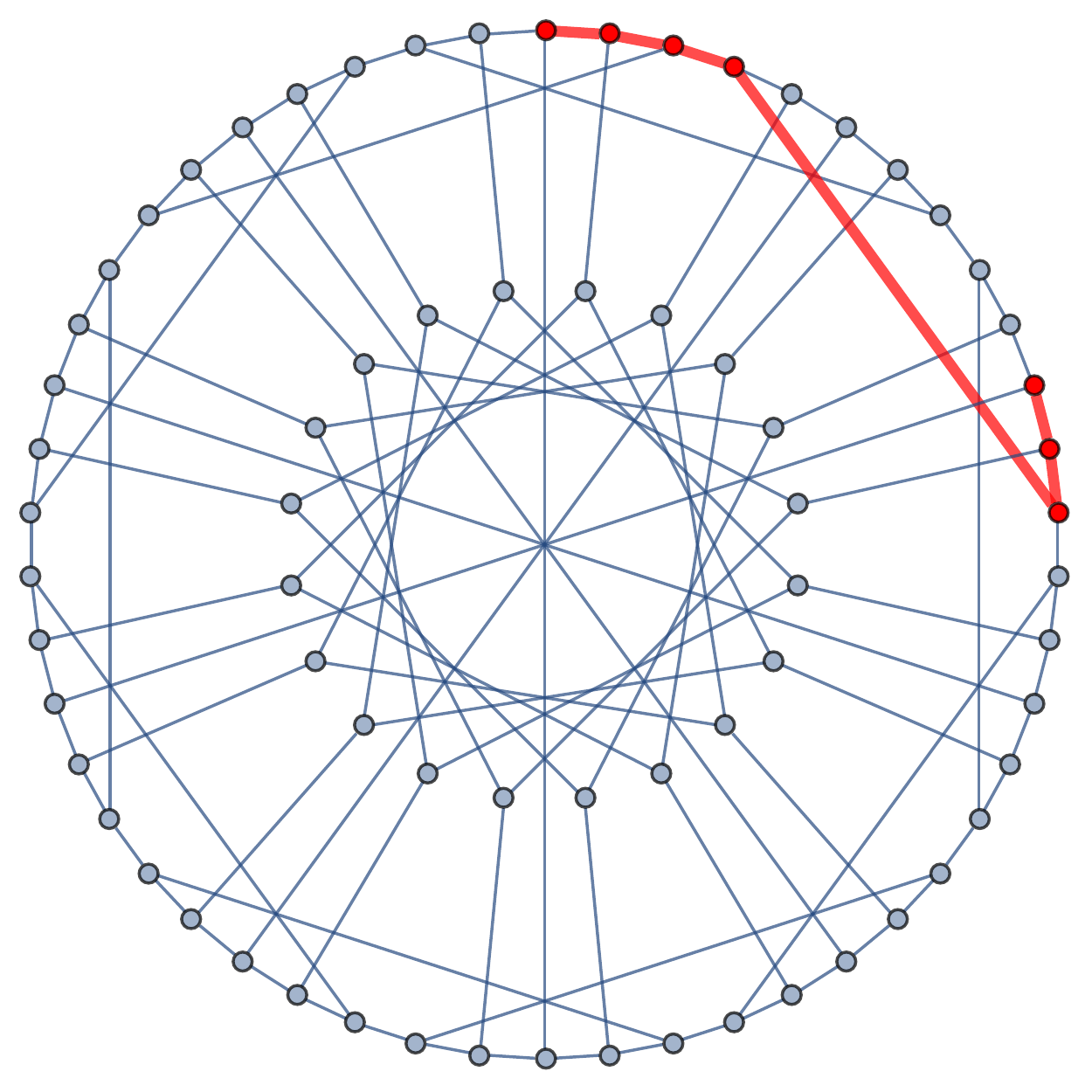}};
\end{tikzpicture}
\end{minipage} 
\caption{Paths taken by the Spectral Method.} 
\end{figure}

\subsection{Related Results and Problems.} We first discuss some explicitly related approaches in Computer Graphics and then relate these questions to classical problem for Partial Differential Equations.\\

\textbf{Geodesics in Heat.} To the best of our knowledge, the philosophically most related approach that we were able to find in the literature is the heat construction of Crane, Weischedel \& Wardetzky \cite{crane, crane1} to obtain a natural notion of distance. Their idea consists of solving the heat equation with an initial point source $\delta_x$ for a short amount of time $t$ to obtain $e^{t\Delta} \delta_x$. Varadhan's Lemma \cite{var, var2} suggests that this solution behaves roughly like a Gaussian in the geodesic distance. This information is then used by solving a second elliptic equation to obtain a normalized distance. The method is shown to be effective in obtaining a notion of intrinsic distance on a variety of domains and has inspired a lot of subsequent work, we mention 
\cite{bel, gao, linb, litman, sharp, tao, yang, yang2}. Crane, Weischedel \& Wardetzky \cite{crane, crane1} also establish that the method can, in principle, recover the exact distances on graphs as $t \rightarrow 0$ (the implicit constant in their result was shown by the author to have combinatorial significance \cite{stein}). This may not be reasonable in practice since an accurate computation of $t \rightarrow 0$ requires the use of all eigenvectors and very high accuracy. Other PDE-type methods have been used, most notably the eikonal equation $|\nabla u|=1$ \cite{cam, kimmel, memoli, memoli2, sethian},
the Schrodinger equation \cite{rang}, the biharmonic equation \cite{lip} and smoothed distances using spectral geometry \cite{coifman}.
Our method is quite different since we do not work with the heat equation at all and instead used a pinned ground state of the Laplacian. In particular, our method is not directly inspired by a continuous method (where imposing that the solution of an elliptic equation vanishes in a single point comes with some difficulty in $d=2$ and significant issues stemming from potential-theoretic obstructions when $d \geq 3$). \\

\textbf{Hot Spots Conjecture.} The Hot Spots conjecture, first posed by Rauch at a 1974 conference in Tulane, asks whether, for a domain $\Omega \subset \mathbb{R}^2$, the first nontrivial eigenfunction of the Laplacian with Neumann boundary conditions assumes maximum and minimum on the boundary. It is not quite clear what to assume on the domain $\Omega$.
 Burdzy \& Werner \cite{b4} constructed a counterexample for a domain that is not simply connected, so the conjecture is now commonly stated either for simply connected domains or convex domains. The Hot Spots conjecture has inspired a lot of work, we refer to \cite{b0, atar, b3, ban, bass, bell, b25, b4, freitas,grieser, grieser2, jerison, jerison2, judge, kawohl, melas, miya, miya2, pascu, manas, steini} and references therein. It is expected that some form of the statement is also true in higher dimensions and on manifolds.
 The relationship to the problem we consider is as follows: let $-\Delta \phi_2 = \mu_2 \phi_2$ denote the second Neumann eigenfunction. It has mean value 0 and, generically, we expect the nodal line $\left\{x: \phi_2(x) = 0\right\}$ to have one connected component that touches the boundary $\partial \Omega$ in two points and to separate $\Omega$ into two connected components $\Omega = \Omega_1 \cup \Omega_2$. This is Payne's nodal domain conjecture \cite{payne}
 which is known to hold in convex domains \cite{al, melas} and known to fail in some other settings \cite{freitas, hoff}. Let us now take a simply connected component of $\left\{x: \phi_2(x) > 0\right\}$. This domain is partially bounded by Dirichlet conditions $\phi_2(x) = 0$ and, on the boundary of $\partial \Omega$, by Neumann boundary conditions $\partial u/\partial n = 0$. This is the natural continuous analogue of the equation we consider graphs. What one observes in the continuous setting is that $\phi_2$ is a remarkably simple object that seems to, more or less, nicely increase with distance from the nodal set; the maxima being on the boundary being merely one instance of this phenomenon. This phenomenon is poorly understood.\\

 One of the contributions of our paper is a simple algorithm on graphs whose efficiency is rooted in the very same phenomenon. In particular, it seems that the underlying reason why the Hot Spots Conjecture is true, is robust enough to survive the transition to graphs. This is well in line with the assumption that versions of the Hot Spots conjecture will also hold in higher dimensions. We point out that Hot Spots problems are not usually considered on graphs due to the difficulty of defining the boundary; we believe our approach to be a suitably adapted analogue that avoids this problem altogether. It is quite conceivable that the discrete setting, allowing for a different set of tools, may give rise to new ideas for the continuous problem and this is part of the motivation underlying this paper. While the existence of the phenomenon, the solution of the PDE behaving nicely with respect to the distance from the nodal set, is not surprising, perhaps its \textit{strength}, mirrored in these surprisingly accurate results on topologically highly nontrivial graphs, is.
 
 \subsection{Other Laplacians} There are various notions of a Laplacian on a graph. Throughout this paper, we will only be concerned with the Laplacian matrix $L=D-A$. However, there are at least two other canonical notions: the random walk normalized Laplacian given by $\mbox{Id}_{n \times n}-D^{-1}A$ and the symmetric normalized Laplacian $\mbox{Id}_{n \times n} - D^{-1/2}AD^{-1/2}$. It is not difficult to see that our proof of Theorem 1 has a natural analogue for these notions and it is possible to use the Algorithm with these other eigenfunctions as underlying geometry. We note that one has to be careful in using the symmetric normalized Laplacian: evaluating the corresponding eigenfunction $\phi$ of $\mbox{Id}_{n \times n} - D^{-1/2}AD^{-1/2}$ at a vertex $v$ results in
 $$ \phi(v) - \sum_{(v,w) \in E}{ \frac{\phi(w)}{\sqrt{\mbox{deg}(v)\mbox{deg}(w)}}} = \mu_2 \phi(v) > 0.$$
 Multiplying the equation with $\mbox{deg}(v)^{-1/2}$ results in
  $$ \frac{\phi(v)}{\sqrt{\mbox{deg}(v)}} - \frac{1}{\mbox{deg}(v)}\sum_{(v,w) \in E}{ \frac{\phi(w)}{\sqrt{\mbox{deg}(w)}}} > 0$$
which shows that every vertex $v \neq i$ has a neighboring vertex $w$ such that 
$$\phi(w) \mbox{deg}(w)^{-1/2} < \phi(v) \mbox{deg}(v)^{-1/2}$$
 and the algorithm has to account for this change and minimize this expression instead. What we observe in basic numerical experiments is that in many cases that we considered, these three Laplacians perform in a very similar way but that $D-A$ seems, generically, to yield slightly better results. It is not clear to us whether there are other differential operators that might end up yielding even better results.
 
 \subsection{Open Problems}
We conclude with a list of open questions; we will only ask them in the discrete setting but they are (with the proper modifications) also interesting in the continuous setting.

\begin{enumerate}
\item \textbf{Rigorous results.} Which rigorous results can be proven about the algorithm? We emphasize once more that any such results might have interesting counterparts in the continuous setting.
\item \textbf{Upper Bounds on the Length.} Is it possible to obtain upper bounds on the length of such an arising path? This would be of interest both in the discrete case but also, for gradient flows, in the continuous setting where it seems to be related to upper bounds on the size of level sets: this questions is notoriously difficult: the main result in this direction were given by Dong \cite{dong}, Donnelly \& Fefferman \cite{don, don1}, Hardt \& Simon \cite{hardt} and Lin \cite{lin}. A natural conjecture in convex domains is that the length of the gradient flow is at most $ c\cdot \mbox{diam}(\Omega)$. It is less clear what to expect on graphs.
\item \textbf{Average Guarantees.} Even though the Algorithm can fail, it seems to fail by very little -- moreover,  in the rare case when it fails by a lot, this seems to occur for a very small number of vertices (see \S 3). Is it possible to obtain upper bounds for a `typical' pair of points?
\item \textbf{Fast Computation.} How fast can this algorithm be carried out? In practice, we are not interested in an arbitrary eigenfunction but the smallest one. Is it important that it is exactly the smallest eigenfunction or are, in practice, small errors permissible? 
\item \textbf{Symmetrization.} The algorithm is not symmetric: the path it finds from $i$ to $j$ will not necessarily be the reverse path that it produces from $j$ to $i$. This allows us to take the better of the two for a more stable estimate. 
\item \textbf{Effective Variants.} Are there effective variations on this idea?  
\end{enumerate}

\begin{center}
\begin{figure}[h!]
\begin{tikzpicture}[scale=0.87]
\node at (0,0) {\includegraphics[width=0.25\textwidth]{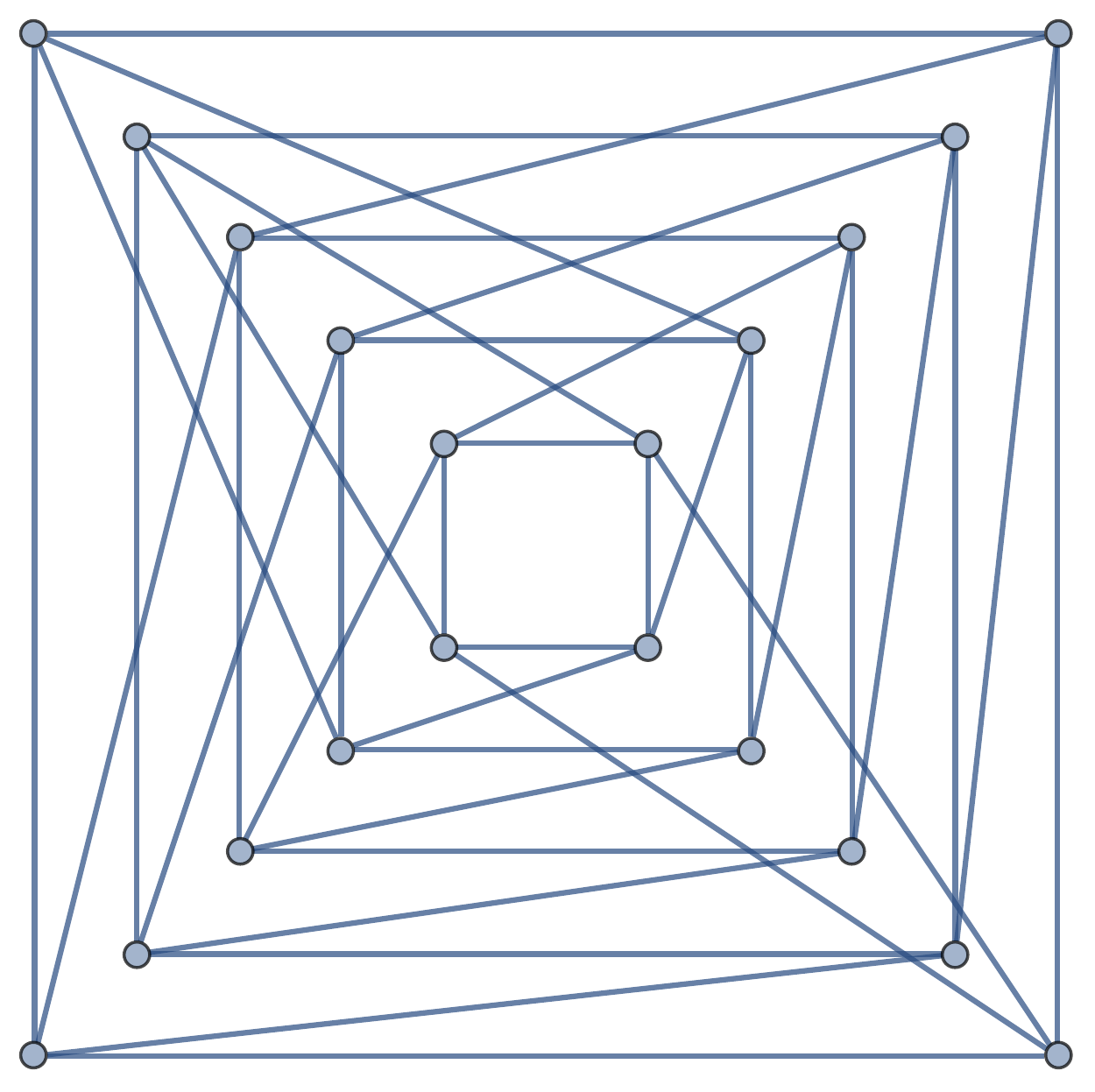}};
\node at (4,0) {\includegraphics[width=0.25\textwidth]{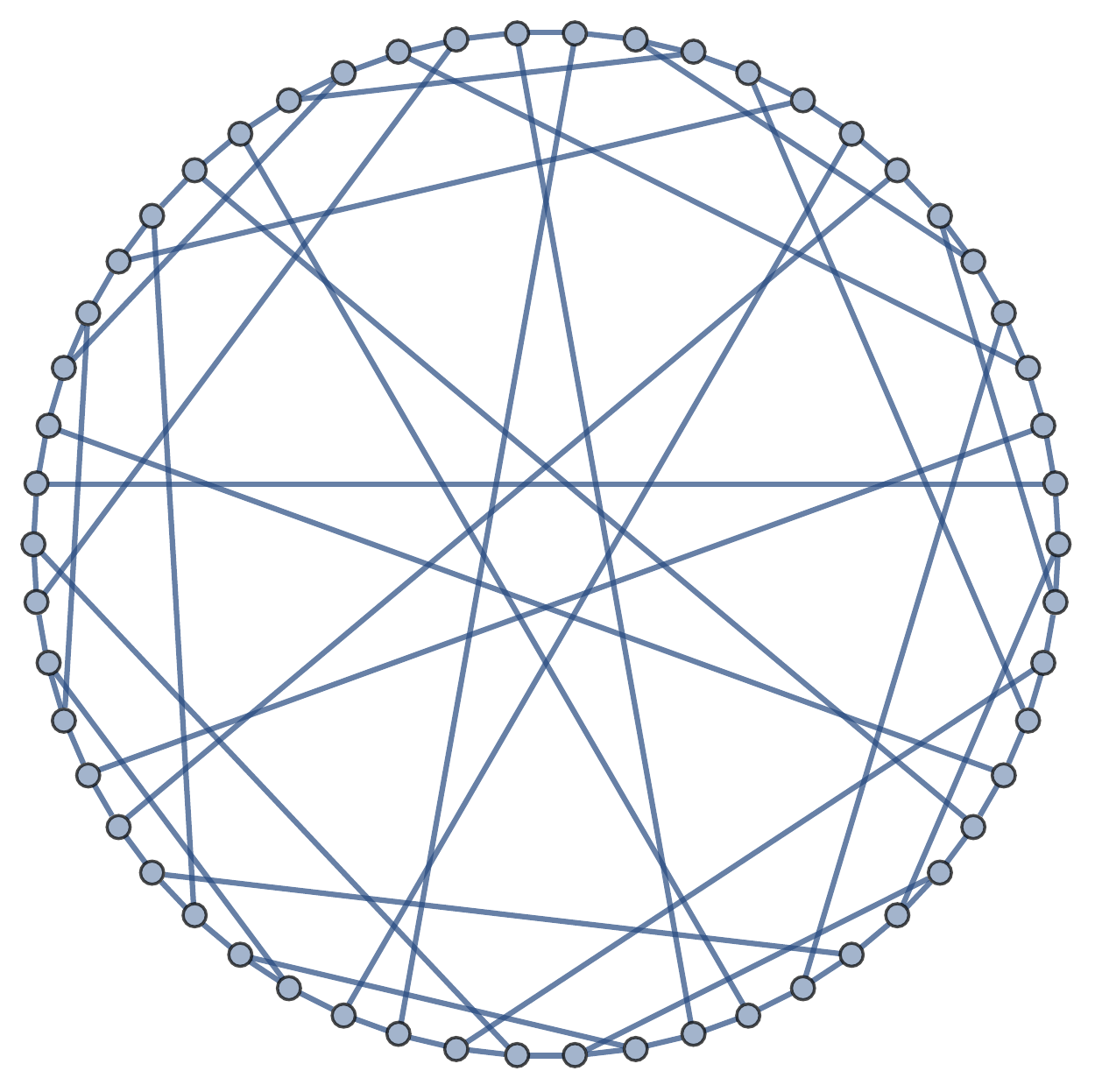}};
\node at (8,0) {\includegraphics[width=0.25\textwidth]{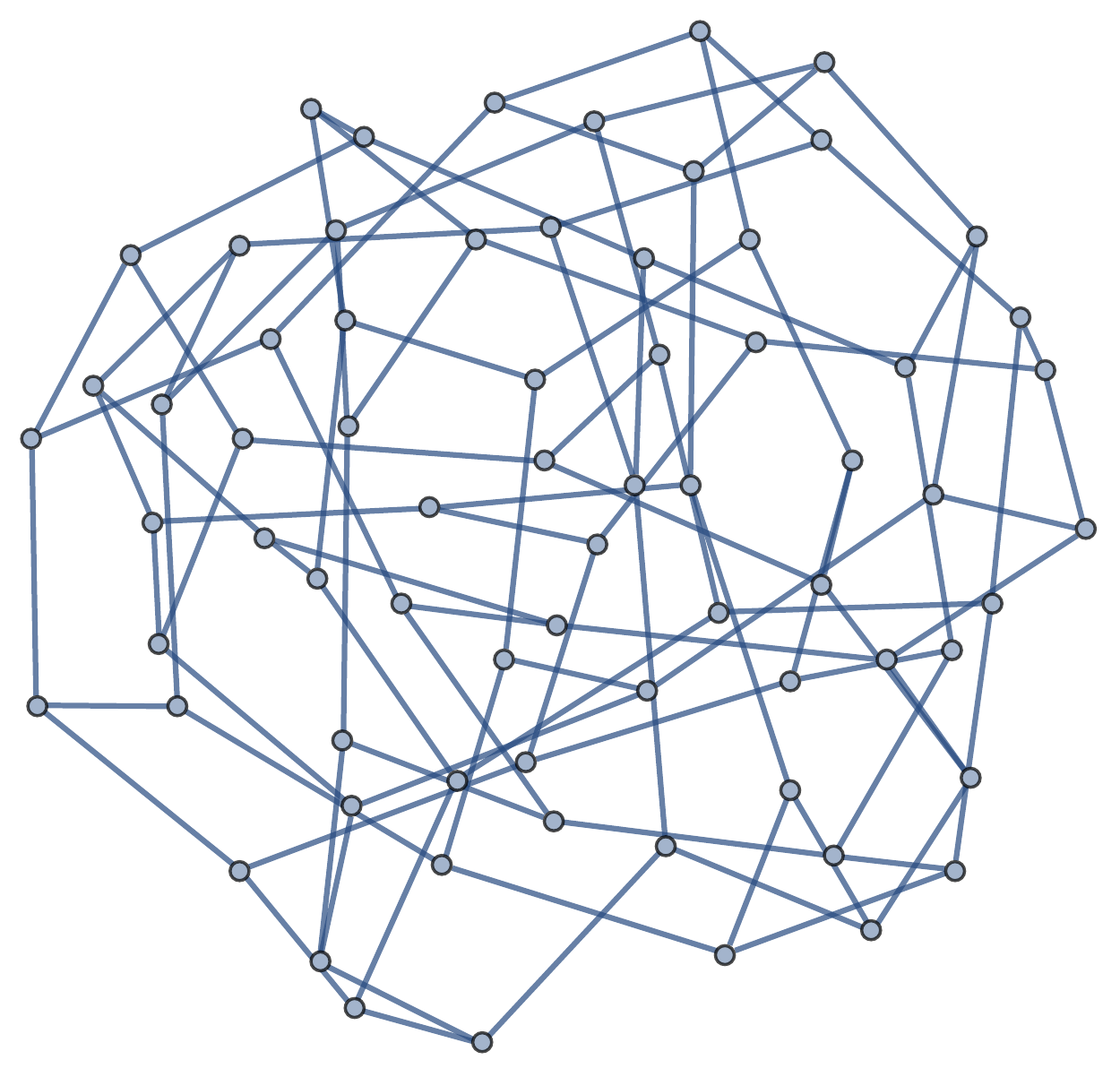}};
\node at (0,-4) {\includegraphics[width=0.25\textwidth]{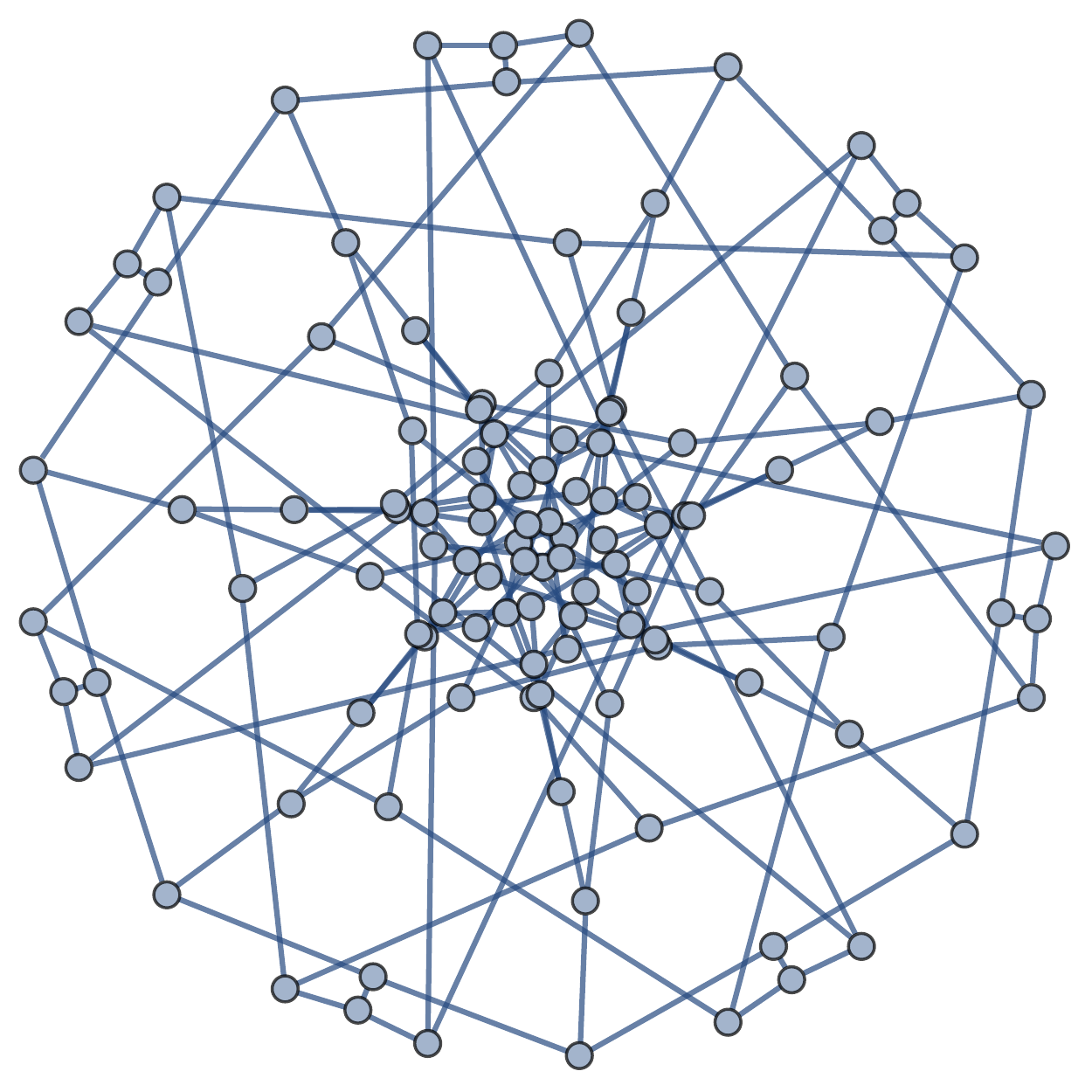}};
\node at (4,-4) {\includegraphics[width=0.25\textwidth]{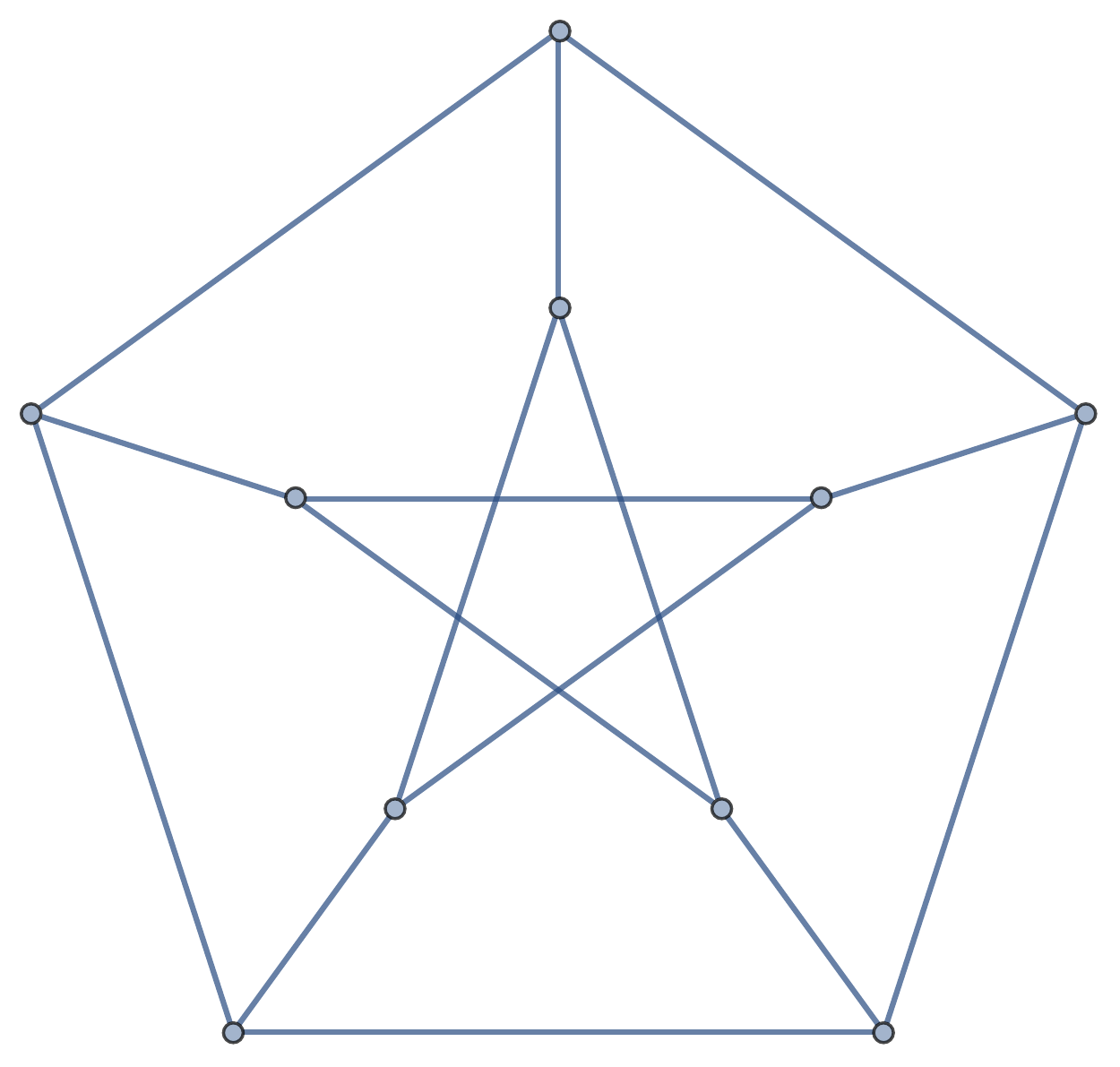}};
\node at (8,-4) {\includegraphics[width=0.25\textwidth]{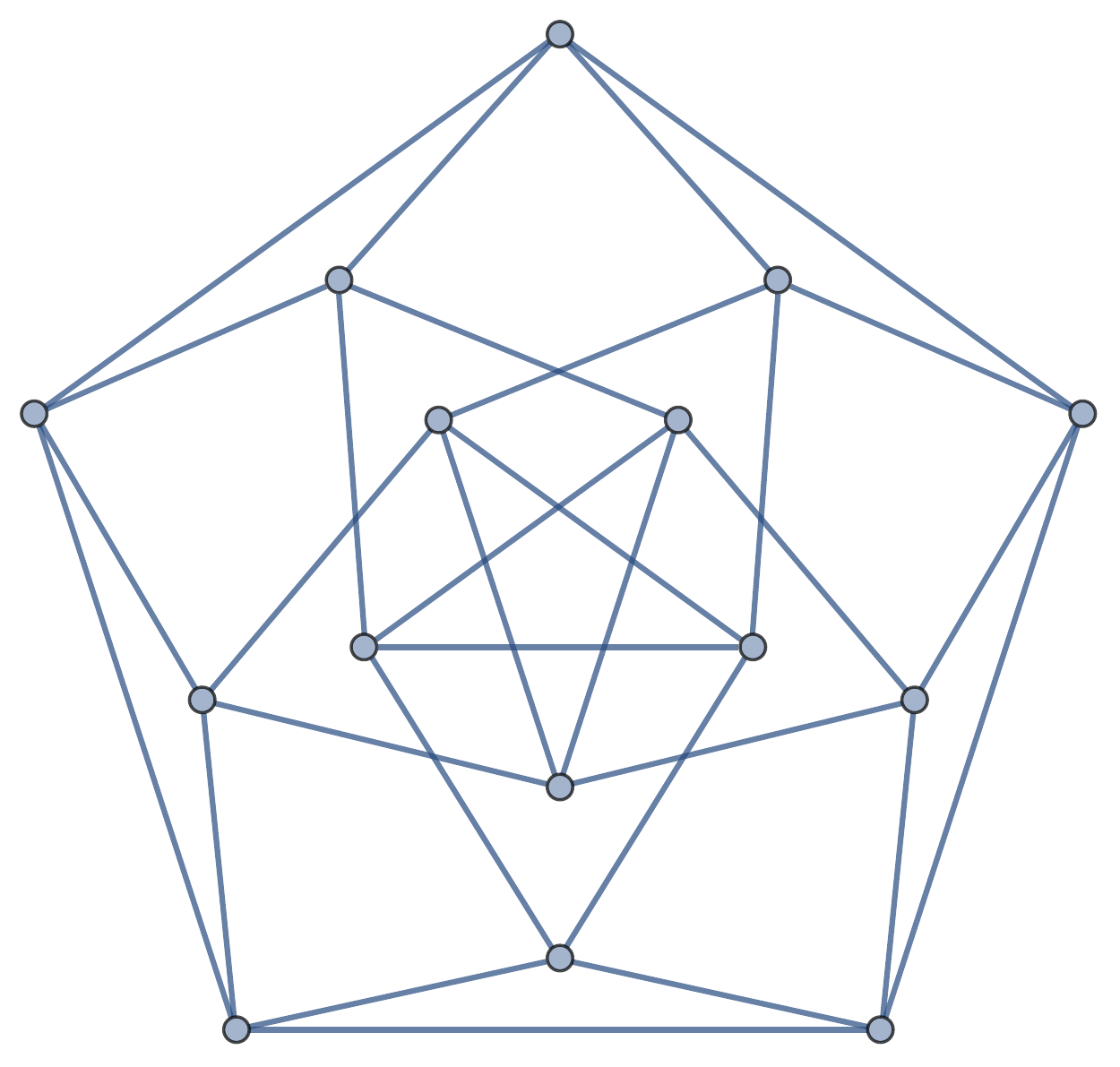}};
\node at (0,-8) {\includegraphics[width=0.25\textwidth]{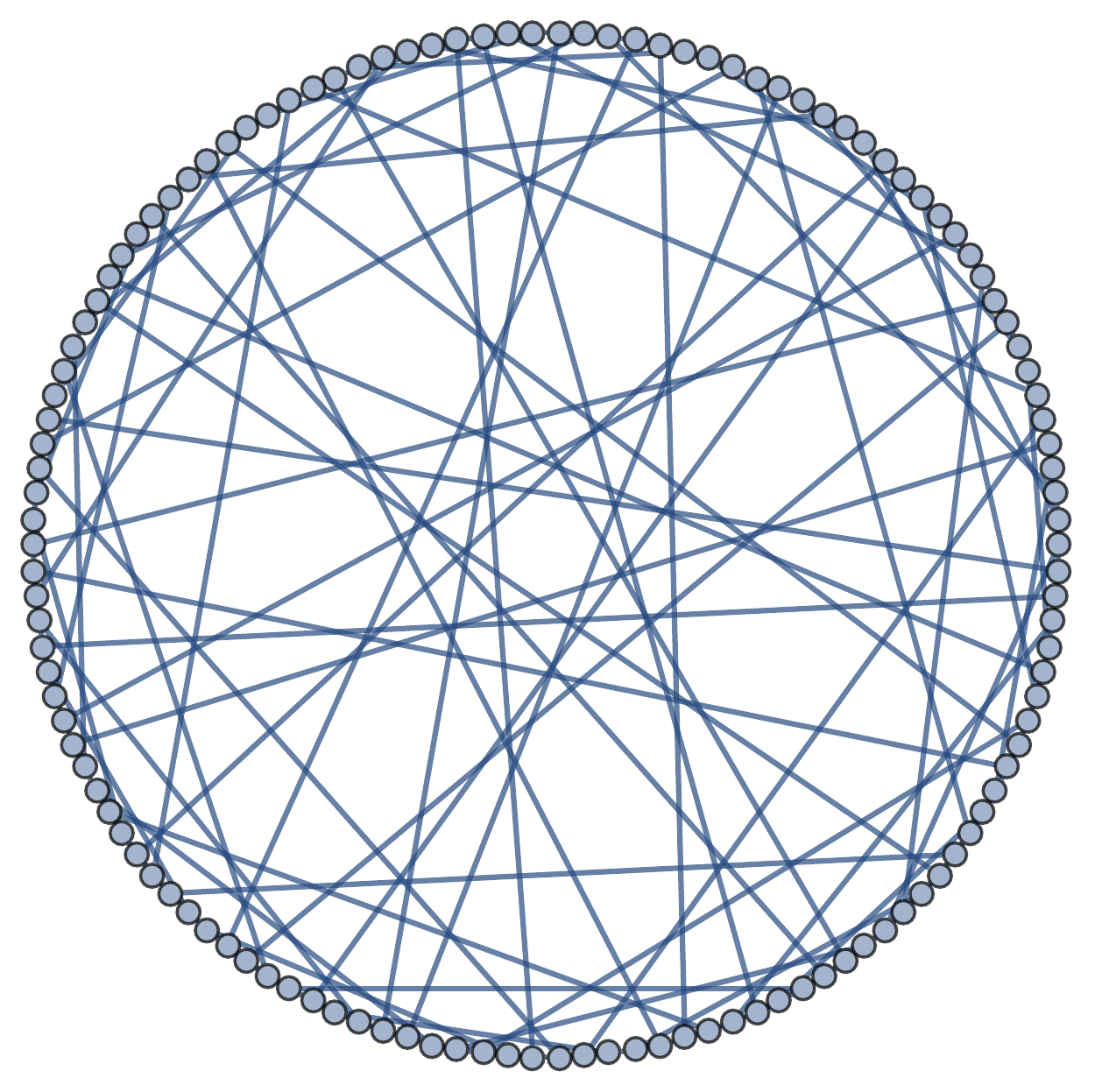}};
\node at (4,-8) {\includegraphics[width=0.25\textwidth]{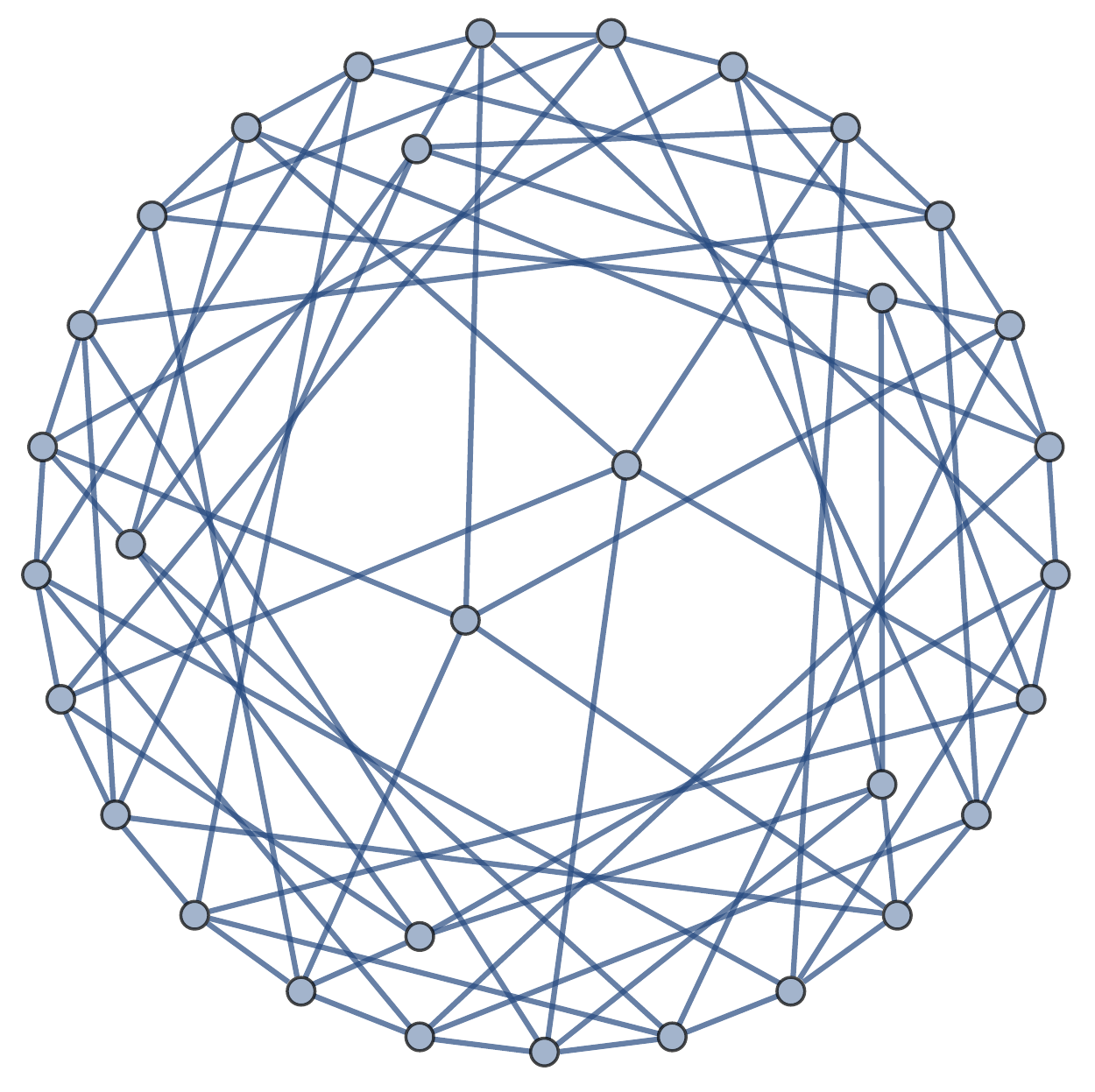}};
\node at (8,-8) {\includegraphics[width=0.25\textwidth]{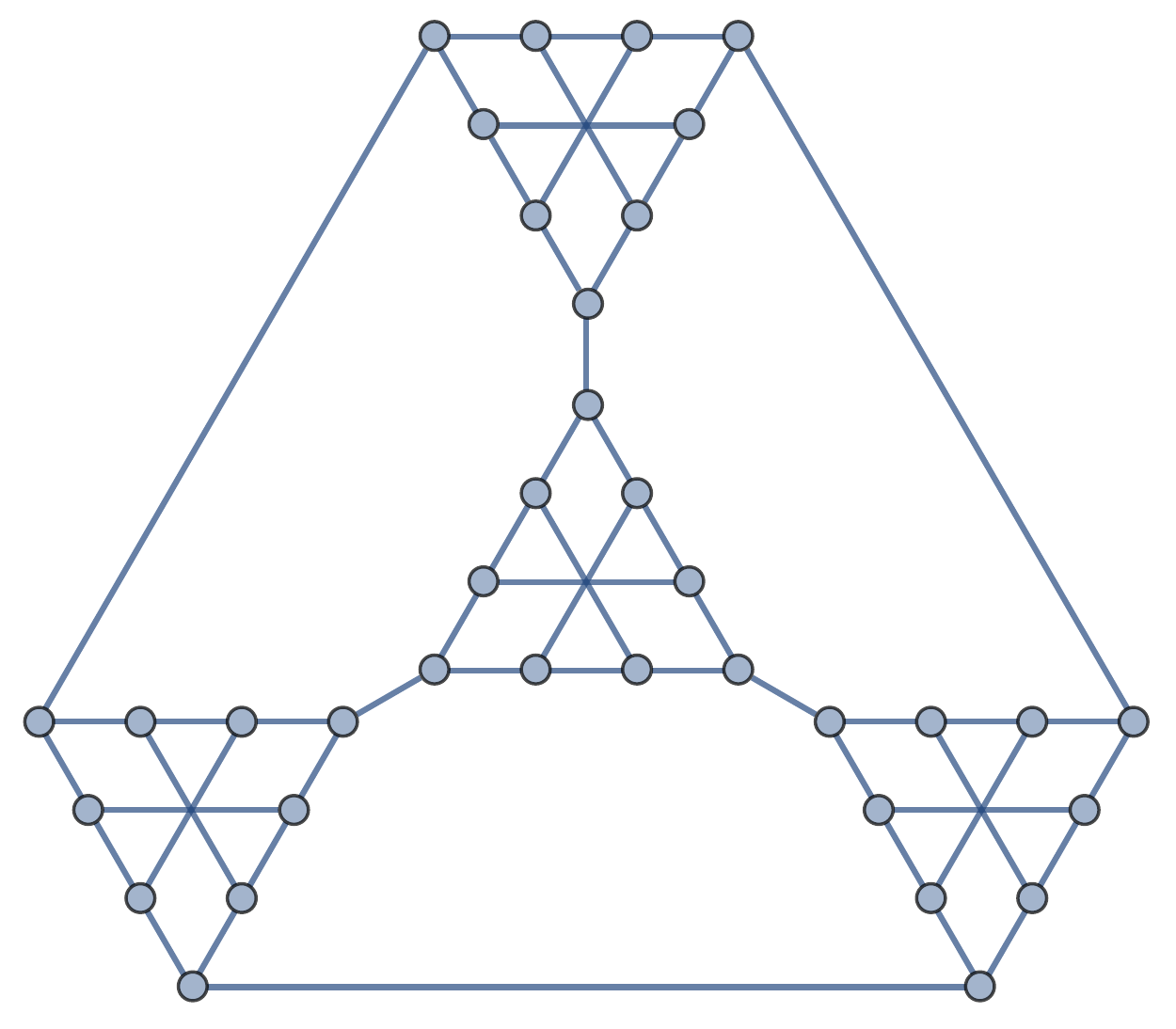}};
\end{tikzpicture}
\caption{Very different types of graphs on which all pairwise shortest paths are reconstructed exactly by the spectral method. The graphs are (1) the Folkman Graph, (2) the Gray graph, (3) the Harries-Wong Graph, (4) the Ljubljana Graph, (5) the Petersen Graph, (6) the Line Petersen Graph, (7) the Tutte 12 Cage, (8) the Wells Graph and (9) the Zamfirescu Graph 36.}
\end{figure}
\end{center}

\section{Some Examples}
We describe some explicit examples. Our observations can be roughly summarized as follows: (1) on graphs with symmetries or some overall homogeneity (see Fig. 4 for some examples), the algorithm tends to produce the shortest path for all pairs of vertices, (2) on graphs that have different densities and bottlenecks, the algorithm still produces the shortest path for most pairs of vertices and very short paths on averages and (3) employing the trick of computing both the optimal path from $i$ to $j$ and from $j$ to $i$ leads to near-optimal results. We see (see \S 3.3 and \S 3.4.) that even in graphs designed to break the algorithm, it yields optimal results for $\sim 95\%$ of all pairs of vertices; moreover, the average path produced by the spectral method is merely $\sim 0.1$ steps longer than the optimal one.

\subsection{Some Examples that work.} We emphasize that the algorithm does, on average, perform extraordinarily well. Even when it seems to fail, it is often on a few pairs of vertices by a little bit. Before discussing some of its failures, we give a list of graphs where it perfectly reconstructs all pairwise distances in Fig. 4.

\subsection{Erd\H{o}s-Renyi Graph.} We continue with the Erd\H{o}s-Renyi model $G(n,p)$. We pick $n$ vertices and connect each of them with likelihood $p$. Numerical Experiments for $100 \leq n \leq 500$ and $p \geq 0.1$ seem to indicate that the method performs perfectly accurately (with high likelihood). Indeed, in this regime we did not find a single example where it fails (though it does fail by small amounts for small $p$). This could be an interesting starting point for a theoretical investigation.

\subsection{Deterministic Graphs} One of the more extreme examples of failure of the algorithm that we found occurs on an amusing graph. On April 1 (!), 1975, a year before the four color theorem was proven, Martin Gardner announced in his column in the \textit{Scientific American} that the `most sensational of last year's discoveries in pure mathematics was surely the finding of a counterexample to the notorious four-color-map conjecture' \cite{gardner} and provided such a map that is now known as Gardner's Map or the April Fools' map. It contains a pair of vertices such that the path arising from the Spectral Method requires 22 steps, 5 steps more than the maximum. The structure of this example is most intriguing, see Figure 5.
\begin{figure}[h!]
\begin{tikzpicture}
\node at (0,0) {\includegraphics[width   = 0.5\textwidth]{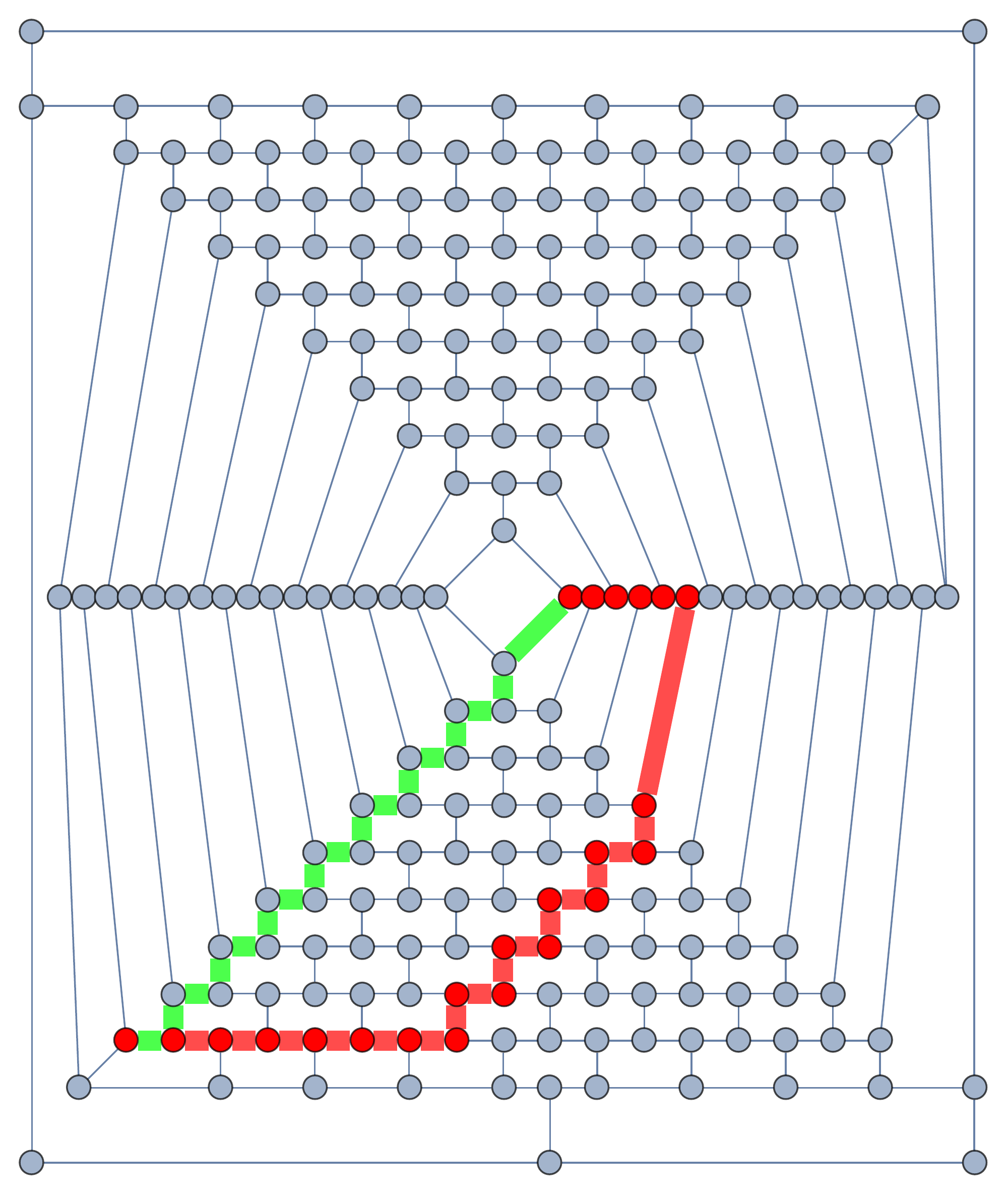}};
\end{tikzpicture}
\caption{An extreme example we found: the shortest path (shown in green) has length 17, the Spectral Method leads to a path (shown in red) of length 22. Nonetheless, for a \textit{randomly chosen} pair of vertices in the graph, the spectral method produces a path that is only 0.08 steps longer than the optimal length.} 
\end{figure}

Even in this remarkable example, the \textit{average} case is very good: the path obtained by the spectral method for a \textit{random} pair of vertices is only 0.08 steps longer than the shortest path.
 We give a second example that occurs on another extreme graph: the Barnette-Bosak-Lederberg graph \cite{lederb} which is the smallest known 3-regular graph without a Hamiltonian cycle, see Fig. 6.
The picture that emerges from these experiments is that counterexamples have to have rather peculiar structure and even then, the failure happens to be concentrated on relatively few points.

 \begin{figure}[h!]
\begin{tikzpicture}
\node at (0,0) {\includegraphics[width = 0.4\textwidth]{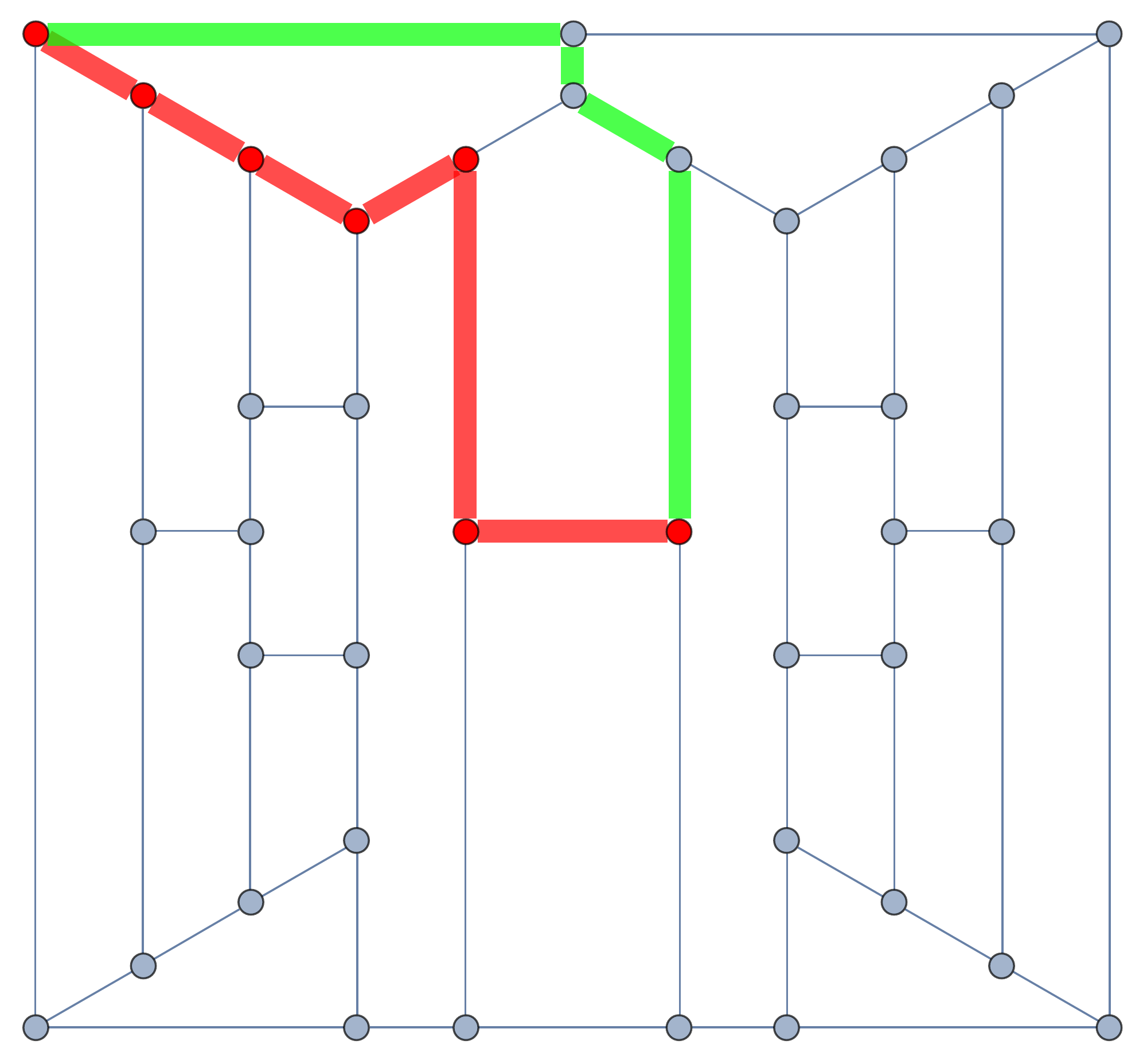}};
\end{tikzpicture}
\caption{The  Barnette-Bosak-Lederberg graph: the shortest path (shown in green) has length 4, the Spectral Method leads to a path (shown in red) of length 6. Nonetheless, for a random pair of vertices in the graph, the spectral method produces a path that is only 0.09 steps longer than the optimal length.} 
\end{figure} We pose two questions (that are also of interest for non-planar graphs):
\begin{enumerate}
\item \textbf{Question 1.} Are there planar graphs such that, for some pair of vertices $u \neq v$, the path produced by the spectral method is $c-$times as long as the shortest path for any $c>1$?
\item \textbf{Question 2.}  Is there a universal constant $c_0$ such that for all planar graphs the path produced by the spectral method between a \textit{random} pair of vertices is at most $c_0-$times as long as the shortest path? 
\end{enumerate}

\begin{figure}[h!]
\begin{tikzpicture}
\node at (0,0) {\includegraphics[width = \textwidth]{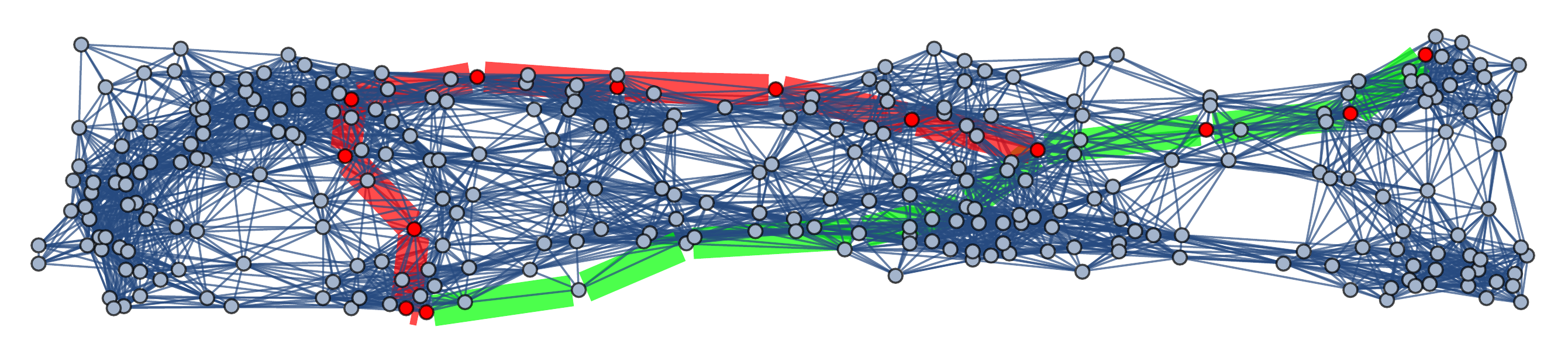}};
\end{tikzpicture}
\caption{The shortest path (shown in green) has length 8, the Spectral Method leads to a path (shown in red) of length 12. For a random pair of vertices in the graph, the spectral method produces a path that is only 0.05 steps longer than the optimal length.} 
\end{figure}

\subsection{Random Geometric Graph.} Consider a random geometric graph: we pick 300 points uniformly at random in the rectangle $[0,3] \times [0,1]$ and connect any two points with an edge if they are distance at most $0.3$. 
 We find that the Spectral Method (run twice: to find a path from $i$ to $j$ and a path from $j$ to $i$) performs exceedingly well on average: the \textit{average} path between two vertices as found by the spectral method is 0.05 steps longer than the shortest path. The most extreme case among all pairs of vertices is shown in the figure: the Spectral Method finds a path of length 12 while the optimal path has length 8.

\section{Proofs}

\subsection{Proof of Theorem 1}
\begin{proof} Proving feasability of the algorithm requires us to show several statements.
\begin{enumerate}
\item The smallest eigenvector $\phi$ of $L_i$ does not change sign (validity of step 2) and can thus be made positive (validity of step 3).
\item If $\phi(v) \neq 0$, then there exists a neighbor $(w,v) \in E$ of $v$ with $\phi(w) < \phi(v)$. 
\item $\phi$ only vanishes in the vertex $i$.
\end{enumerate}
We start by analyzing the variational problem
$$ \phi = \arg \min_{f(i) = 0} \sum_{(w_1, w_2) \in E}{(f(w_1)-f(w_2))^2} \qquad \mbox{subject to}~ \sum_{v \in V}{f(v)^2} = 1.$$

Clearly, a minimizer exists by compactness. We observe two cases: either the graph becomes disconnected when removing $i$ or it stays connected. If it becomes disconnected (which, for example, would be the case for trees), then we can think of the problem as being comprised of two or more graphs being connected via $i$. It thus suffices to prove the statement for each component individually and we can assume that $G$ remains connected after removing $i$. We shall assume this henceforth.
Since
$$ \sum_{(w_1, w_2) \in E}{(f(w_1)-f(w_2))^2} = \left\langle f, L_i f\right\rangle$$
and $L_i$ is a symmetric matrix, we know that a minimizer corresponds to the smallest eigenvalue of $L_i$. This eigenvalue is easily seen to be bigger than 0: if it were 0, connectedness of the graph implies that $f$ is constant
which creates a contradiction between its $\ell^2-$norm being 1 and $f(i)=0$. Let us now suppose that $f$ is such a minimizer. 
We define $g(\cdot) = |f(\cdot)|$ and note that $g(i) = 0$, the sum of the squares of $g$ is still 1 and
\begin{align*}
 \sum_{(w_1, w_2) \in E}{(g(w_1)-g(w_2))^2} &= \sum_{(w_1, w_2) \in E}{(|f(w_1)|-|f(w_2)|)^2}\\
 &\leq  \sum_{(w_1, w_2) \in E}{(f(w_1)-f(w_2))^2}.
 \end{align*}
However, since $f$ was a minimizer, we have to have equality: for all $(w_1, w_2) \in E$ we have to have
$$ \left| f(w_1)-f(w_2)\right| = | |f(w_1)| - |f(w_2)||.$$
This, in turn, requires that if $(w_1, w_2) \in E$ are two connected vertices, then $f(w_1)$ and $f(w_2)$ have the same sign (though it is conceivable that one of the values is 0). In particular, if $f$ were to change sign, then, since the Graph remains connected after removing $i$, there would have to be at least one vertex $k$ where $f$ and thus $g$ vanishes. In summary, we have obtained a function $g$ such that (1)
$$  \sum_{(w_1, w_2) \in E}{(g(w_1)-g(w_2))^2} =  \sum_{(w_1, w_2) \in E}{(f(w_1)-f(w_2))^2}$$
is minimal, (2) $g \geq 0$ and (3) $g$ vanishes in at least one additional vertex that is not $i$. Moreover, again by connectedness, there exists one vertex $\ell$ such that $g$ vanishes in $\ell$ and $\ell$ has at least one neighbor where $g$ does not vanish. We conclude with one elementary fact about real numbers: for any $\left\{x_1, \dots, x_n\right\} \subset \mathbb{R}$ and $x \in \mathbb{R}$, we have
$$ \sum_{k=1}^{n}{ (x_k - x)^2} \geq  \sum_{k=1}^{n}{ \left(x_k - \frac{1}{n}\sum_{i=1}^{n}{x_i}\right)^2}$$
with equality if and only if $x$ is indeed the mean of the $n$ numbers. This means that if we replace the value of $g$ in $\ell$ by the mean of the neighbors, then the Dirichlet energy decreases and the $L^2-$norm increases which is a contradiction to $g$ minimizing the functional. Therefore $f$ does not change sign and only vanishes in $i$. This implies uniqueness: suppose there was a second minimizer, $g$, distinct from $f$. Since both $f$ and $g$ can be interpreted as eigenvectors of a symmetric matrix, we know that they can be chosen to be orthogonal. However, since a minimizer does not change sign, the inner product of any two minimizers cannot vanish which establishes uniqueness.
 It remains to show that the algorithm cannot get `stuck'. Let $v \neq i$ be a vertex and suppose the algorithm finds itself in $v$. We take the $v-$th row of the equation $L_i f = \mu f$ and obtain that
$$ 0 < \mu f(v) = \sum_{(v,w) \in E}{ (f(v) - f(w))} =\mbox{deg}(v)f(v) -  \sum_{(v,w) \in E}{ f(w)}.$$
Therefore
$$ \sum_{(v,w) \in E}{ f(w)} < \mbox{deg}(v)f(v)$$
and there exists at least one neighbor of $v$ where $f$ assumes a smaller value.
This establishes the desired result. \end{proof}

\begin{center}
\begin{figure}[h!]
\begin{tikzpicture}[scale=0.6]
\filldraw (-3,0) circle (0.05cm);
\filldraw (-3,1) circle (0.05cm);
\filldraw (-2,0) circle (0.05cm);
\filldraw (-1,0) circle (0.05cm);
\filldraw (0,0) circle (0.05cm);
\filldraw (1,0) circle (0.05cm);
\filldraw (1,1) circle (0.05cm);
\filldraw (1,-1) circle (0.05cm);
\draw [thick] (-3,0) -- (1,0);
\draw [thick] (1,-1) -- (1,1);
\draw [thick] (-2,0) -- (-3,1);
\draw [thick, dashed] (-3,1) -- (-4,1);
\draw [thick, dashed] (-4,0) -- (-3,0);
\draw [thick, dashed] (1,1) -- (2,1);
\draw [thick, dashed] (1,-1) -- (2,-1);
\node at (0, -0.4) {$i$};
\node at (1.5,0) {$G_2$};
\node at (-3, -0.5) {$G_1$};
\draw[ultra thick] (0,-3) circle (0.2cm);
\filldraw (0,-3) circle (0.05cm);
\filldraw (1,-3) circle (0.05cm);
\filldraw (1,-3) circle (0.05cm);
\filldraw (1,-3) circle (0.05cm);
\draw [thick] (-1,-3) -- (1,-3);
\draw [thick] (1,-4) -- (1,-3);
\draw [thick] (1,-3) -- (1,-2);
\filldraw (1,-2) circle (0.05cm);
\filldraw (1,-4) circle (0.05cm);
\draw [thick, dashed] (1,-2) -- (2,-2);
\draw [thick, dashed] (1,-4) -- (2,-4);
\draw [thick, dashed] (-1,-2) -- (-2,-2);
\draw [thick, dashed] (-1,-4) -- (-2,-4);
\draw [thick] (-1,-2) -- (-1,-4);
\filldraw (-1,-4) circle (0.05cm);
\filldraw (-1,-2) circle (0.05cm);
\node at (0, -3.4) {$i$};
\node at (4.5,0) {};
\end{tikzpicture}
\caption{The Duplication Procedure from the Proof of Theorem 2.}
\end{figure}
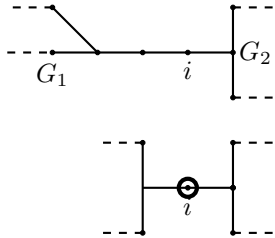
\end{center}

\subsection{Proof of Theorem 2}
\begin{proof} The proof is based on a duplication trick. Let $i \in V$ be given and let
$$ \phi = \arg\min_{f:V \rightarrow \mathbb{R} \atop  f(u) = 0, f \not\equiv 0} \frac{\sum_{(w_1, w_2) \in E}{(f(w_1)-f(w_2))^2}}{\sum_{w \in V}{f(w)^2}}.$$
We know from the proof of Theorem 1 that $\phi$ does not change sign and can w.l.o.g. be taken positive. Since the graph is a tree, we know that removal of the vertex $i$ leads to two connected graphs $G_1$ and $G_2$ (one of them might be empty if $i$ is a leaf). It remains to prove the desired result for $G_1$ and $G_2$ separately. We instead argue with a separation trick (see Fig. 8): we mirror one of the graphs and define its values at a mirrored vertex to be the negative value of $\phi$ at the original vertex and call the arising function $\psi$. 
We see that $\psi$ satisfies $(D-A)\psi = \mu_2 \psi$ since the equation is valid on each side of the graph and is easily verified in $i$ since $\psi(i) = 0$. By fairly standard arguments, we obtain that $\psi$ is the smallest nontrivial eigenfunction on the new graph and we then obtained the desired property from a result of Lefevre \cite{lef} (see als Gernandt \& Pade \cite{gern}). An alternative proof could be given as follows: we can interpret the function $\phi$ on $G_2 \cup \left\{i \right\}$ as a
function satisfying $(D-A) \phi = \mu_2 \phi$ in all vertices except $i$. This allows for $\phi$ to be reinterpreted as the expected outcome of a game that was described by Lederman and the author in \cite{leder}. The game is obviously strictly monotonic in the distance from $i$ and this yields the desired result. We refer to \cite{leder} for the details.
\end{proof}

\begin{center}
\begin{figure}[h!]
\begin{tikzpicture}
\draw [ultra thick] (0,0) -- (0,1);
\draw [ultra thick] (-0.02,0) -- (-0.02,1);
\draw [very thick] (0,0) to[out=0, in=120] (1,-1) to[out=300, in=180] (2, -2) to[out=0, in=300] (1.5, -1) to[out=120, in=180] (1.6, 0) to[out=0, in=110] (2.5, -1) to[out=290, in =180] (3.5, -1.5) to [out=0, in =270] (4, -1) to[out=90, in=0] (3.5, 0) to[out=180, in = 180] (3.5, -0.5) to[out=0, in=90] (3.6, -1) to[out=270, in=300] (3, -1) to[out=150, in =180] (4, 1) to[out=0, in=0] (2, 2) to[out=180, in=45] (0,1);
\node at (-0.5, 0.5) {$\partial {\Omega_0}$};
\node at (3.2, 2) {$\partial {\Omega_n}$};
\filldraw (2,1) circle (0.05cm);
\node at (2.2, 0.8) {$x_0$};
\draw [thick, ->] (2, 1) -- (2.2, 1.2);
\draw [thick, ->] (2, 1) -- (2-0.2, 1-0.2);
\draw [thick, dashed] (2.2, 1.2) to[out=45, in=235] (3.3, 1.5);
\end{tikzpicture}
\caption{Sketch of the Proof of Theorem 3}
\end{figure}
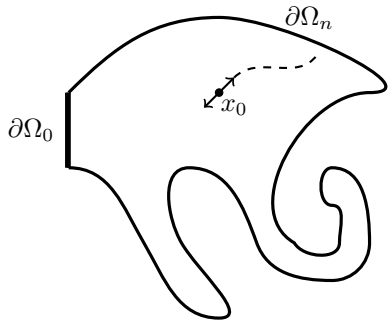
\end{center}
\subsection{Proof of Theorem 3.}
\begin{proof} It follows from standard arguments that $u > 0$ inside the domain. Let us now suppose that there is a critical point $x_0$ inside the domain. Since we are in two dimensions and $\Delta u < 0$, the critical point has to be either a local maximum (in which case it lies on the boundary by assumption) or it is an isolated critical point $x_0$. By the Morse Lemma, the function $u$ has a local saddle structure. In particular, there exist two directions in which $u$ is increasing. We start two increasing gradient (ascent) flows $\gamma_1$ and $\gamma_2$, one in each direction. 
These gradient ascent flows can never intersect themselves or each other (we explain below why). This means they ascend until they hit a critical point. This critical point can be a maximum (in which case it is at the boundary and we stop) or it is another critical point: in that case, we use the local saddle structure, rotate the vector by $90^{\circ}$ and continue the gradient ascent until we eventually hit a maximum on the boundary. These flows can never hit $\partial \Omega_0$ because the function is increasing along them -- this means they cut the domain into two parts, one containing $\partial \Omega_0$, the other only containing $\partial \Omega_n$. We pick the second part, $\Omega_2$, and note that $\partial u/\partial n = 0$ on the entire domain: on $\partial \Omega_n$ this is merely the imposed boundary condition, on the part given by $\gamma_1$ and $\gamma_2$, this follows from the fact that the tangent line of $\gamma_i$ is parallel to the gradient of $u$ in the point. Then, however,
$$  0 = \int_{\partial \Omega_2} \frac{\partial u}{\partial n} d\sigma = \int_{\Omega_2}{ \Delta u} < 0$$
which is a contradiction.  The same argument could have also been applied if the gradient ascent curves were to intersect themselves or each other.
It remains to show that the Gradient flow converges to $\Omega_0$. We pick a point $x_0$ and start the gradient flow $\gamma_1$ and, simultaneously, the gradient ascent flow $\gamma_2$. We know that $\gamma_2$ terminates in a critical point and these are all on the boundary. If $\gamma_1$ were to not terminate in $\partial \Omega_0$ but on
$\partial \Omega_n$, then we could again extract a domain for which $u$ has Neumann boundary conditions and obtain a contradiction. This concludes the argument.
\end{proof}

\end{document}